\documentclass[11pt]{article}

\usepackage{latexsym,mathrsfs}
\usepackage{amsmath,amssymb}
\usepackage{amsthm,enumerate,verbatim}
\usepackage{amsfonts}
\usepackage{graphicx}
\usepackage{algorithm}
\usepackage{algorithmic}
\usepackage{url}

\newtheorem{theorem}{Theorem}
\newtheorem{proposition}{Proposition}
\newtheorem{corollary}{Corollary}
\newtheorem{definition}{Definition}
\newtheorem*{definition*}{Definition}

\newtheorem{remark}{Remark}

\numberwithin{equation}{section}
\numberwithin{table}{section}
\numberwithin{figure}{section}
\def \R{{\mathbb R}}
\def \C{{\mathbb C}}

\newcommand {\mat}  [1] {\left[\begin{array}{#1}}
\newcommand {\rix}      {\end{array}\right]}

\DeclareMathOperator{\argmin}{argmin}

\title{Approximating the nearest stable discrete-time system}

\author{
Nicolas Gillis\thanks{Department of Mathematics and Operational Research, University of Mons, Rue de Houdain 9, 7000 Mons, Belgium. Email: \{nicolas.gillis, punit.sharma\}@umons.ac.be. NG and PS acknowledge the support of the ERC (starting grant n$^\text{o}$ 679515). NG also acknowledges the support of the F.R.S.-FNRS (incentive grant for scientific research n$^\text{o}$ F.4501.16).} 
\qquad 
Michael Karow\thanks{Institut f${\rm \ddot{u}}$r Mathematik, MA 4-5 TU Berlin, Str.\@ d.\@ 17.\@ Juni 136, D-10623 Berlin, Germany.
Email: karow@math.tu-berlin.de.} 
\qquad 
Punit Sharma\thanks{Department of Mathematics, Indian Institute of Technology Delhi, Hauz Khas, New Delhi-110016, India; \texttt{punit.sharma@maths.iitd.ac.in}.
P. Sharma acknowledges the support of the DST-Inspire Faculty Award (MI01807-G) by Government of India and Institute SEED Grant (NPN5R) by IIT Delhi.}
}

\begin{document}

\maketitle

\begin{abstract}
In this paper, we consider the problem of stabilizing discrete-time linear systems by computing a nearby stable matrix to an unstable one. To do so, we provide a new characterization for the set of stable matrices. We show that a matrix $A$ is stable if and only if it can be written as $A=S^{-1}UBS$, where $S$ is positive definite, $U$ is orthogonal, and $B$ is a positive semidefinite contraction (that is, the singular values of $B$ are less or equal to 1). This characterization results in an equivalent non-convex  optimization problem with a feasible set on which it is easy to project. We propose
a very efficient fast projected gradient method
to tackle the problem in variables $(S,U,B)$ and generate locally optimal solutions.
We show the effectiveness of the proposed method compared to other approaches.
\end{abstract}


\section{Introduction}
Consider a discrete-time linear system described by the following difference equation
\begin{equation}\label{dessys}
x(t+1)=Ax(t),\quad t \in \mathbb N,
\end{equation}
where $A\in \R^{n,n}$ and $\mathbb N$ is the set of nonnegative integers,
$x(t)$ denotes the $n$-dimensional state vector.
If $\lambda_1,\ldots,\lambda_n$ are the eigenvalues of $A$, then such a system
is called stable (resp.\@ asymptotically stable) if $|\lambda_i|\leq 1$ (resp.\@ $|\lambda_i|<1$) for all $i=1,\ldots,n$, and the eigenvalues
with unit modulus are semisimple; otherwise, it is called unstable.\\
In this paper, we consider the \emph{nearest stable matrix problem} in the discrete-time case.
More precisely, for a given unstable matrix $A$, we consider the following optimization problem
\begin{equation}\label{eq:probdef}
\inf_{X\in \mathbb S_d^{n,n}}{\|A-X\|}_{F}^2,
\end{equation}
where ${\|\cdot\|}_F$ denotes the Frobenius norm of a matrix and
$\mathbb S_d^{n,n}$ is the set of all stable matrices of size $n \times n$.
We consider in this paper the Frobenius norm of the error as it is arguably the most widely used norm.
However, our approach can be directly applied to any differentiable cost function
(e.g., any component-wise $\ell_p$ norm with $p > 1$).

\paragraph*{Notation} Throughout the paper, $X^T$, ${\rm tr}(X)$, and
$\|X\|$  stand for the transpose, the trace and the spectral norm of a real square matrix $X$, respectively.
By $\Lambda(X)$, $\rho(X):=\max_{\lambda \in \Lambda(X)}{|\lambda|}$ and $\kappa(X)=\|X\|\, \|X^{-1}\|$ we denote the spectrum (set of eigenvalues), the spectral radius and the condition number. We write
 $X\succ 0$ and $X\succeq 0$ $(X \preceq 0)$ if $X$ is symmetric and positive definite
or positive semidefinite (symmetric negative semidefinite), respectively. The positive semidefinite symmetric square root of a positive semidefinite symmetric matrix $X$ is denoted by $X^{1/2}$.

\paragraph*{Related work} For a given unstable matrix $A$ (in the discrete- or continuous-time case), the problem
of computing the smallest perturbation that stabilizes $A$, also known as
 the nearest stable matrix  problem, occurs in system identification
where one needs to identify a stable system depending on observations~\cite{OrbNV13}.
To the best of our knowledge, the nearest stable matrix problem was first introduced and analyzed
in the discrete- and continuous-time case in~\cite{OrbNV13}, where a nearby stable approximation $X$ of a given unstable system $A$ is constructed by means of successive convex approximations of the set of stable systems.
For the continuous-time case, two methods were recently proposed:
\begin{enumerate}

\item In~\cite{GilS17}, the problem is reformulated into an equivalent problem with a simple convex feasible set. In fact, it is shown that $A$ is stable  if and only if it can be written as $A=(J-R)Q$ where $J^T=-J$, $R\succeq 0$, and $Q \succ 0$.
This result was later generalized to solve various nearness problems for continuous-time linear systems, namely, the problems of finding
the nearest stable Metzler matrix~\cite{anderson2017distance},
the nearest stable matrix pair~\cite{GilMS17}
and the nearest positive real system~\cite{GilS17b}.

\item  In~\cite{GugL17}, the problem is tackled by solving low-rank matrix differential equations.

\end{enumerate}

The nearest stable matrix problem in discrete-time case has received much less attention, and to the best our knowledge, only~\cite{OrbNV13} considered this problem without any assumption on the entries of the matrix.
For the class of positive systems of the form~\eqref{dessys}, where the matrix $A$ is component-wise nonnegative, the problem of computing the nearest stable nonnegative matrix has been studied very recently in~\cite{NesP17,GP2018}. In~\cite{NesP17}, authors consider the nearest stable/unstable nonnegative matrix with respect
to the max-norm ${\|X\|}_{\max} = \max_{i,j}|X_{i,j}|$,
the $\ell_{\infty}$ operator norm ${\|X\|}_{\infty} = \sup_{u \neq 0} \frac{{\|Xu\|}_{\infty}}{{\|u\|}_{\infty}}$, and
 the $\ell_{1}$ operator norm ${\|X\|}_{1} = \sup_{u \neq 0} \frac{{\|Xu\|}_{1}}{{\|u\|}_{1}}$,
where ${\|x\|}_{\infty} = \max_i |x_i|$ and ${\|x\|}_{1} = \sum_i |x_i|$.
For these norms, it turns out that, rather surprisingly, the problem can be solved in polynomial-time.
In~\cite{GP2018}, authors propose a more efficient heuristic than in~\cite{OrbNV13}
for the Frobenius norm for which the problem is more difficult, with the existence of many local minima (up to $2^n$ in dimension $n$). \\
The nearest stable matrix problem~\eqref{eq:probdef} is the converse
problem of stability radius problem in the discrete-time case, where a stable matrix $A$ is given and one looks for the smallest perturbation that moves an eigenvalue outside the stability region. The converse problem has been introduced and studied extensively; see, e.g., \cite{GenAV02,HinP90,HinS98,HinSN03,NgoN03,NgoN05} and the references therein. \\
The problem~\eqref{eq:probdef} is notoriously difficult due to properties of the spectral radius as a function
of matrix: the set $\mathbb S_d^{n,n}$ of stable matrices is highly nonconvex~\cite{OrbNV13}, and neither open
nor closed. For example, $B_\epsilon \notin \mathbb S_d^{2,2}$ for $\epsilon >0$ but $B \in \mathbb S_d^{2,2}$, where
{\small
\[
\underbrace{\mat{cc}
1 & \epsilon\\
-\epsilon & 1 \\
\rix}_{=:B_\epsilon}
\rightarrow
\underbrace{\mat{ccccc}
1 & 0 \\
0 & 1 \\
\rix}_{=:B},
\]}
and
$C_\delta \in \mathbb S_d^{2,2} $ for $0\leq \delta < 1$, but $C \notin \mathbb S_d^{2,2}$, where
{\small
\[
\underbrace{\mat{cc}
1 & 1\\
0 & \delta\\
\rix}_{=:C_\delta} \rightarrow
\underbrace{\mat{ccccc}
1 & 1 \\
0 & 1 \\
\rix}_{=:C}.
\]
}
Therefore it is in general difficult to obtain a global optimal solution to problem~\eqref{eq:probdef}.

The aim of this paper is to derive counterparts of a number of results in~\cite{GilS17} for the discrete-time
case. These results require special constructions and show special features in the discrete-time case.
Our principle strategy for computing a nearby stable approximation to a given unstable matrix is to
reformulate the problem~\eqref{eq:probdef} into an equivalent problem with a simple feasible set onto which
points can be projected relatively easily. We aim to provide in many cases a better approximation than the one obtained
with the code from~\cite{OrbNV13} at a lower computational cost.

The paper is organized as follows. In Section 2, we define the SUB form of a matrix: the matrix $A$ admits a SUB form if there exists $(S,U,B)$ with $S\succ 0$, $U$ is orthogonal, $B\succeq 0$ and $\|B\|\leq 1$ such that $A = S^{-1}  U B S$.
We prove that a matrix is stable if and only if it admits a SUB form.
In Section 3, we propose a fast gradient method (FGM) to solve the
reformulated problem in variables $(S,U,B)$, along with several initialization strategies.
To illustrate the performance of FGM, we apply it on several examples of unstable matrices and compare the results with the algorithm from~\cite{OrbNV13}.

\section {A new characterization for stable matrices} \label{sec:reformulation}

In this section, we derive a factorization of stable matrices into symmetric and orthogonal matrices.
This will allow us to reformulate problem~\eqref{eq:probdef} into an equivalent problem with a
simple feasible set for which standard optimization methods can be applied. In order to do this,
we define the SUB form of a matrix.
\begin{definition}
A matrix $A \in \mathbb R^{n,n}$ is said to admit a \emph{SUB form} if there exist $S,U,B \in \mathbb R^{n,n}$
such that $A=S^{-1}UBS$ where $S\succ 0$, $U$ is orthogonal, $B\succeq 0$ and $\|B\|\leq 1$.
\end{definition}
\begin{theorem}\label{thm:mainresult}
A matrix is stable (asymptotically stable) if and only if it admits a SUB form (SUB form with $\|B\| < 1$).
\end{theorem}
%
\par\noindent{\bf Proof}. \ignorespaces
The proof follows by the following two facts: 1) the Lyapunov criterion of the Schur stability~\cite{Gan59a}; 2) the polar decomposition.
Indeed, by the Lyapunov theorem, $A$ is stable (asymptotically stable) if and only if there exists
an ellipsoid $E$ such that $AE \subset E$ (respectively, $AE \subset \text{int}E$). This is equivalent to say that there
exist matrices $C$ and $L$ such that $\|L\| \leq 1$ (respectively, $\|L\| < 1$) and $A = C^{-1}LC$. Now we write
the polar decomposition $C = V S$, where $V$ is orthogonal and $S \succ 0$. Thus, $A = S^{-1}V^{-1}LV S$.
Denote $V^{-1}LV = M$. Clearly, $\|M\| = \|L\|$. Finally, write the polar decomposition: $M = UB$
with $U$ orthogonal, $B \succ 0$, and $\|B\| = \|M\|$. We have $A = S^{-1}UBS$, which completes the proof.
\space{\ \vbox{\hrule\hbox{\vrule height1.3ex\hskip0.8ex\vrule}\hrule}}\par
In view of Theorem~\ref{thm:mainresult}, the set $\mathbb S_d^{n,n}$ of stable matrices can be
characterized as the set of matrices that admit a SUB form, or equivalently, we can parameterize the set of stable matrices using a matrix triple $(S,U,B)$ as follows
\begin{eqnarray*}
\mathbb S_d^{n,n}=\Big\{S^{-1}UBS \in \R^{n,n}~\big|~S\succ 0,~U~\text{orthogonal},
~B\succeq 0~\text{with}~
\|B\|\leq 1
\Big\}.
\end{eqnarray*}
This characterization changes the feasible set and the objective function in the
nearest stable matrix problem~\eqref{eq:probdef} as
\begin{eqnarray}\label{eq:thm_reform}
\inf_{X\in \mathbb S_d^{n,n}}{\|A-X\|}_{F}^2 =
\inf_{S\succ 0,~U\,\text{orthogonal},~B \succeq 0,~ \|B\| \leq 1}{\|A-S^{-1}UBS\|}_{F}^2.
\end{eqnarray}
As we mentioned earlier, the set $\mathbb S_d^{n,n}$ of stable matrices
is neither open nor closed and clearly the new parameterization of $\mathbb S_d^{n,n}$
in terms of matrix triple $(S,U,B)$ does not change this, since $\mathbb S_d^{n,n}$
is not open because of the constraint $B\succeq 0$ and not closed due to the constraint
$S \succ 0$.
Therefore the infimum in the right hand side of~\eqref{eq:thm_reform} may not be attained.
%

In the next section, we will provide an algorithmic solution for the
nearest stable matrix problem~\eqref{eq:probdef} by trying to solve
the reformulated problem~\eqref{eq:thm_reform}.
%
%
Note that in view of Lyapunov's Theorem replacing $S$, $U$, and $B$ in terms of the variable $P\succ 0$ leads to the
formulation
\[
\inf_{X,\,P\succ 0} {\|A-X\|}_F\quad \text{such that}\quad X^T P X - P \preceq 0,
\]
which is difficult to solve numerically as it involves highly non-linear constraints~\cite{OrbNV13}. Thus a key contribution
of this paper is the reformulation~\eqref{eq:thm_reform}.
An advantage of this reformulation is that the feasible set is rather simple and
therefore it is relatively easy to project onto it. As a result we propose
a fast projected gradient method to solve the reformulated problem~\eqref{eq:thm_reform}, see Algorithm~\ref{fastgrad}.
We close the section with some useful remarks on the matrices that admit a SUB form.
\begin{remark}\label{rem:alhpa-stab}{\rm
Let $\gamma \in \R$ with $0<\gamma < 1$. Then $A$ is called $\gamma$-stable if $\lambda \in \Lambda(A)$
satisfies $|\lambda|\leq \gamma$. Note that $A$ is $\gamma$-stable if and only if $B=\frac{A}{\gamma}$
is stable, since for any nonzero $x \in \mathbb C^n$ we have $Ax=\mu x$ if and only if  $\mu=\gamma \lambda$ for some
$\lambda \in \Lambda(B)$ such that $Bx=\lambda x$.
Thus
from Theorem~\ref{thm:mainresult} $A$ is $\gamma$-stable if and only if $\frac{A}{\gamma}$ admits a SUB form
if and only if $A$ admits a SUB form with $\|B\|\leq \gamma$. This observation can be used to find
a nearby $\gamma$-stable matrix to a given unstable one.
}
\end{remark}

\begin{remark}\label{rem:uniqueness}{\rm
We note that a remark for the non-uniqueness of the SUB decomposition of a stable matrix similar to
\cite[Remark 6]{GilS17} for continuous systems also holds in the discrete case. The SUB representation
of a stable matrix $A$, that is, $A=S^{-1}UBS$ where $S\succ 0$, $U^TU=I_n$, $B\succeq 0$, and $\|B\| \leq 1$,
is non-unique. As there is always a scaling degree of freedom: for any scalar $\alpha > 0$ we have
$A=(\alpha S)^{-1}UB(\alpha S)$. This can partially be addressed by the fact that the ellipsoid $E$ in the Lyapunov Theorem is non-unique and the matrices $S$ and $B$ in the SUB form
depend on $E$, see the proof of Theorem~\ref{thm:mainresult}. However, characterizing precisely the
non-uniqueness of the SUB form (and possibly taking advantage of it in a numerical algorithm) is
a direction for further research.
}
\end{remark}

\begin{remark}\label{rem:complex}{\rm
The results of this section are readily extended to handle complex matrices by
substituting $X^*$, the conjugate transpose for $X^T$ and unitary matrices for
orthogonal matrices. In particular, we have that
$A\in \C^{n,n}$ is stable if and only if there exist $S,U,B \in \mathbb C^{n,n}$
such that $A=S^{-1}UBS$ where $S\succ 0$, $U$ is unitary, $B\succeq 0$ and $\|B\|\leq 1$.
We note  that a similar observation also holds for the characterization of complex stable matrices in the
continuous-time case. In particular following the terminology in~\cite{GilS17} we have that $A\in \C^{n,n}$
is stable in the continuous-time case if and only if there exist $J,R,Q\in \C^{n,n}$ such that
$A=(J-R)Q$ where $J^*=-J$, $R\succeq 0$, and $Q \succ 0$. This was not mentioned in~\cite{GilS17}.
}
\end{remark}

\section{Algorithmic solutions to the nearest stable matrix problem}\label{sec:algo}

As shown in Section~\ref{sec:reformulation}, finding the nearest stable matrix to an unstable one is equivalent to solving~\eqref{eq:thm_reform}.
%
In this section, we propose a fast projected gradient method~\cite[p.90]{Nes04} to tackle~\eqref{eq:thm_reform}.
Although fast gradient methods (FGM's) were initially designed for convex optimization problems,
they have recently been shown to work well for non-convex problems as well; see, e.g., \cite{GL16,agarwal2017finding,o2017behavior}.
In particular, for the problem of finding the nearest stable matrix in the continuous-time case,
they work significantly better than standard gradient schemes and coordinate descent methods~\cite{GilS17}, while being relatively simple to implement.
We use a similar implementation as in~\cite{GilS17b}; see Algorithm~\ref{fastgrad} for the details.
As for the standard projected gradient method, FGM requires the computation of the gradient of the objective function, and the projection onto the feasible set.
The gradient of $f(S,U,B) = {\|A-S^{-1}UBS\|}_{F}^2$ with respect to $S$ is given by
\[
\nabla_S f(S,U,B) = 2 \, S^{-T} [R^T (R-A)  - (R-A)R^{T}],
\]
where $R = S^{-1}UBS$.
The details are given in~\ref{appgradS}.
For $U$ and $B$, we have
\[
\nabla_U f(S,U,B) = - 2 S^{-1} (A-R) S B^T
\]
and
\[
\nabla_B f(S,U,B) = - 2 U^T S^{-1} (A-R) S .
\]
The projections of a solution $(S,U,B)$ onto the feasible set of~\eqref{eq:thm_reform} are described in Section~\ref{projsec}.

\algsetup{indent=2em}
\begin{algorithm}[ht!]
\caption{Fast Gradient Method (FGM) for \eqref{eq:thm_reform} with restart from~\cite{GilS17b}} \label{fastgrad}
\begin{algorithmic}[1]
\REQUIRE
An initialization $X = (S,U,B)$,
a parameter $\alpha_1 \in (0,1)$,
a lower bound for the step length $\underline{\gamma}$,
an initial step length $\gamma > \underline{\gamma}$.

\ENSURE An approximate solution $X = (S,U,B)$ to~\eqref{eq:thm_reform}.  \medskip

%

\STATE $X' = X$. \emph{\% Create the second sequence of iterates of FGM.}

\FOR{$k = 1, 2, \dots$}

\STATE  $\hat{X} = X$. \emph{\% Keep the previous iterate in memory.}

\STATE \emph{\% Project the gradient step from $X'$ onto the feasible set.}

\STATE $X = \mathcal{P}(X' - \gamma \nabla f(X') )$. \emph{\%  $\mathcal{P}$ is the projection onto the feasible set.}

\STATE    \emph{\% Check if the objective function $f$ has decreased, otherwise decrease the step length.}

\WHILE { $f(X)  > f(\hat{X})$ and $ \gamma \geq \underline{\gamma}$ }

\STATE $\gamma = \frac{2}{3} \gamma$.

\STATE $X = \mathcal{P}( X' - \gamma \nabla f(X') )$.

\ENDWHILE

\STATE  \emph{\% If the step length has reached the lower bound ($f$ could not be decreased), reinitialize $X'$ (the next step will be a standard gradient descent step).}

\IF { $\gamma < \underline{\gamma}$ }

\STATE Restart fast gradient: $X'$ = $X$; $\alpha_{k} = \alpha_{1}$.

\STATE   Reinitialize $\gamma$ at the last value for which it allowed decrease of $f$.

\ELSE

\STATE $\alpha_{k+1} = \frac{1}{2} \left(  \sqrt{ \alpha_{k}^4 + 4 \alpha_{k}^2 } - \alpha_{k}^2 \right)$, $\beta_k =  \frac{\alpha_{k} (1-\alpha_{k})}{\alpha_{k}^2 + \alpha_{k+1}}$.

\STATE $X' = X + \beta_k \left( X - \hat X \right)$.

\ENDIF

\STATE $\gamma = 2 \gamma$.

\ENDFOR

\end{algorithmic}
\end{algorithm}

\paragraph*{Convergence.} Algorithm~\ref{fastgrad} is guaranteed to decrease the objective function at each step because of the line-search (steps 7-10).
Hence, at every iteration, we have
$\| A-S^{-1}UBS\|_F\leq f_0$
where $f_0$ is the initial objective function value.
Since the objective function is bounded from below by zero,
this implies that the objective function values converge to some value $f^*$.
Moreover, the approximations $S^{-1}UBS$ generated at each step of the algorithm are in a compact set : in fact,
\[
{\| S^{-1}UBS \|}_F - {\|A\|}_F  \leq {\| A  - S^{-1}UBS \|}_F \leq f_0
\quad \Longrightarrow \quad
{\|S^{-1}UBS\|}_F \leq f_0+{\|A\|}_F.
\]
Therefore, there exists a subsequence of approximations $S^{-1}UBS$ generated by Algorithm~\ref{fastgrad} that converge to some limit point $A_p^*$
with ${||A- A_p^*||}_F = f^*$.
However, it is more difficult to prove convergence of the iterates $(S,U,B)$ as $S$ is not bounded (e.g., if $A = 0$, then $B = 0$ is optimal while $S$ can be any invertible matrix).
It is possible to add an upper bound on the norm of $S$ to guarantee a subsequence of iterates to converge, but we have not observed in practice that this was an issue.
Providing a rigorous proof of convergence of the iterates of Algorithm~\ref{fastgrad} to a stationary point of~\eqref{eq:thm_reform} is a difficult problem which we leave as a question for further research. It has to be noted that only stationary points are fixed point of our method, since this is a projected gradient method.

\paragraph*{Parameters.} Algorithm~\ref{fastgrad} is not too sensitive to the initial step length $\gamma$ as it increases/decreases it to allow the objective function to decrease,
and reinitialize the value to the previous value that allowed decrease when it is restarted (step 14).
We chose the initial step length to be $\gamma = 1/L$
where $L = \left(\frac{\lambda_{\max}(S)}{\lambda_{\min}(S)}\right)^2 = \kappa(S)^2$.
The reason for this choice is that $L$ is the Lipschitz constant of the gradient of $f$ with respect to $B$ so that using the step length $1/L$ would guarantee the decrease of $f$ if
we would only optimize over $B$ as the problem in variable $B$ is convex~\cite{Nes04}.
For $\alpha_1$, we use 0.5 as in~\cite{GilS17b}.
Since~\eqref{eq:thm_reform} is a difficult non-convex optimization problem, any local optimization scheme such as our FGM approach will be sensitive to initialization; this is discussed in Section~\ref{subsec:int}.

\subsection{Projections} \label{projsec}

In this section we derive the relevant formulas to project $S$, $U$ and $B$ onto the feasible set of~\eqref{eq:thm_reform}.

\paragraph*{Projections for $S$ and $B$.}
 In order to calculate the projection of a square matrix onto the set of positive semidefinite contractions, we introduce some notation.
For a  symmetric matrix $H \in \R^{n,n}$ with eigenvalues $\lambda_k$ ($1 \leq k \leq n$) and orthogonal diagonalization $H=V{\rm diag}(\lambda_1,\ldots, \lambda_n)V^{T}$, we set
$f(H)=V{\rm diag}(f(\lambda_1),\ldots, f(\lambda_n))V^{T}$, where $f$ is any complex valued function defined on the spectrum of $H$. The matrix $f(H)$ does not depend on the particular orthogonal matrix $V$ since it is easily verified that $f(H)=q(H)$, where $q$ is any polynomial that maps each $\lambda_k$ to its value $f(\lambda_k)$.  For a general matrix $X \in \R^{n,n}$, we consider functions of its symmetric part, $f^s(X):=f((X+X^T)/2)$. For an interval $[a,b]\subset \R \cup \{\infty\}$ and $\lambda\in \R$ let
$$
p_{a,b}(\lambda):=\max\{a,\min\{b, \lambda\}\}=
\begin{cases}
a &\text{if }\lambda< a,\\
\lambda &\text{if }\lambda \in [a,b],\\
b&\text{if }b< \lambda.
\end{cases}
$$
Then $p_{a,b}(\lambda)$ is the nearest point projection of $\lambda$ onto $[a,b]$, that is,
$|\lambda-p_{a,b}(\lambda)|={\rm argmin}_{h \in [a,b]}|\lambda-h|$.
%
The statement below extends \cite[Lemma 10]{AlMa12} to the case that $X$ is nonsymmetric.
\begin{proposition}
With respect to Frobenius norm
the matrix $p_{a,b}^s(X)$ is the nearest point projection of $X \in \R^{n,n}$ onto the set
${\mathcal I}_{a,b}=\{\, H \in \R^{n,n}\, |\; H=H^T,\; a \, I\preceq H\preceq b\, I\}$, that is,
 $$p_{a,b}^s(X)={\rm argmin}_{H \in{\mathcal I}_{a,b}} {\|X-H\|_F}.$$
\end{proposition}
\par\noindent{\bf Proof}. \ignorespaces
Let $(X+X^T)/2=V{\rm diag}(\lambda_1,\ldots, \lambda_n)V^T$ with orthogonal $V$. Let $H\in {\mathcal I}_{a,b}$, and let
$\tilde H=V^T H V=[\tilde h_{ij}]$. Then $\tilde H\in {\mathcal I}_{a,b}$ and therefore $\tilde h_{ii}
\in[a,b]$ for all $i=1, \ldots , n$. By orthogonality between symmetric and skew symmetric matrices
and the orthogonal invariance of the Frobenius norm we have
\begin{eqnarray}\label{eqtemp1}
{\|X-H\|}_F^2
&=&{\left\|\frac{X-X^T}{2}\right\|}_F^2+{\left\|\frac{X+X^T}{2}-H\right\|}_F^2 \nonumber \\
&=& {\left\|\frac{X-X^T}{2}\right\|}_F^2+{\|{\rm diag}(\lambda_1,\ldots, \lambda_n)-\tilde H\|}_F^2 \nonumber \\
  &= & {\left\|\frac{X-X^T}{2}\right\|}_F^2
  + \sum_i(\lambda_i-\tilde h_{ii})^2+\sum_{i \not =j}\tilde h_{ij}^2.
\end{eqnarray}
The sum is minimized by
$\tilde H={\rm diag}( p_{a,b}(\lambda_1), \ldots , p_{a,b}(\lambda_n))$. Thus,
$H=p_{a,b}^s(X)$.
\space{\ \vbox{\hrule\hbox{\vrule height1.3ex\hskip0.8ex\vrule}\hrule}}\par
Since for a positive semidefinite matrix the inequality $\|B\| \leq \alpha$ is equivalent to $B\preceq \alpha I_n$
we have the corollaries below.
\begin{corollary}
 With respect to Frobenius norm the nearest point projection of $X \in \R^{n,n}$ onto the set of positive semidefinite
contractions is $p_{0,1}^s(X)$, that is,
$$p_{0,1}^s(X)={\rm argmin}_{B \succeq 0, \|B\| \leq 1} {\|X-B\|_F}.$$
\end{corollary}
\begin{corollary}{\rm \cite{Hig88b}}
 With respect to Frobenius norm the nearest point projection of $X \in \R^{n,n}$ onto the cone of $n\times n$ positive semidefinite
matrices is $p_{0,\infty}^s(X)$, that is,
$$p_{0,\infty}^s(X)={\rm argmin}_{B \succeq 0} {\|X-B\|_F}.$$
\end{corollary}

\paragraph*{Projections for $U$.} Before we give the projection onto the set of orthogonal matrices, we provide another closely related projection that will be useful to obtain initializations in Section~\ref{subsec:int}.
Note that these results require the polar decomposition~\cite{HorJ85}.
\begin{proposition}\label{thm:optimUB}
Let $X\in \R^{n,n}$ and let $X=VH$ be the polar decomposition of $X$, where $V\in \R^{n,n}$ is orthogonal and $H \in \R^{n,n}$ satisfies $ H\succeq 0$. Then
$$
{\argmin}_{(U,B), U^TU=I_n, B \succeq 0, \|B\| \leq 1}
{\| X - UB\|}_F^2 =
\left( V, p_{0,1}(H)\right),
$$
\end{proposition}
%
%
\par\noindent{\bf Proof}. \ignorespaces
Let $H=Q\,{\rm diag}(\lambda_1, \ldots, \lambda_n)Q^T$ be a diagonalization of $H$ with orthogonal $Q$.
Let $U,B\in \R^{n,n}$ be such that $U^TU=I_n$ and $B\succeq 0$ with $\|B\|\leq 1$.
Then
\begin{eqnarray}\label{eq:optUB1}
{\|X-UB\|}_F^2&=&{\|VH-UB\|}_F^2={\|H-V^TUB\|}_F^2 \nonumber\\
&=& {\|Q\, {\rm diag}(\lambda_1, \ldots, \lambda_n)Q^T -V^T UB\|}_F^2 \nonumber\\
&=&
{\|{\rm diag}(\lambda_1, \ldots, \lambda_n)-Q^TV^TUBQ\|}_F^2 \nonumber\\
&\geq & \sum_{i}(\lambda_i-p_{-1,1}(\lambda_i))^2\label{eq:esti5}
\\
&= & \sum_{i}(\lambda_i-p_{0,1}(\lambda_i))^2.\nonumber
\end{eqnarray}
The last equation holds since all $\lambda_i$'s are nonnegative.
The inequality (\ref{eq:esti5}) follows from the fact that all diagonal
entries of $Q^TV^TUBQ$ are contained in $[-1,1]$
since $\|Q^TV^TUBQ\|=\|B\|\leq 1$.
Equality holds in (\ref{eq:esti5}) if and only if
$Q^T V^T UBQ={\rm diag}( p_{0,1}(\lambda_1), \ldots , p_{0,1}(\lambda_n))$.
The latter is equivalent to $UB = Vp_{0,1}(H)$.
\space{\ \vbox{\hrule\hbox{\vrule height1.3ex\hskip0.8ex\vrule}\hrule}}\par
%
%
\begin{proposition}
 Denoting $\mathcal{P}_{\bot}(X)$ the projection of $X$ onto the set of $n\times n$ orthogonal matrices, we have
$
\mathcal{P}_{\bot}(X)
= \argmin_{U^TU = I_n} {\|X-U\|}_F
= V$, where $X = V H$ is the polar decomposition of $X$.
\end{proposition}
The proof is analogous to the proofs of the other propositions in this section and therefore omitted.
%
\subsection{Initializations}\label{subsec:int}

In this section, we propose three initializations.

\paragraph{Standard initialization} We use $S = I_n$, for which the optimal values of $U$ and $B$ can be computed using the polar decomposition of $A$,
see Proposition~\ref{thm:optimUB}:
\[
{\rm argmin}_{(U,B), U^TU=I_n, B \succeq 0, \|B\| \leq 1}
{\| A - UB\|}_F =
\left( V, p_{0,1}(H)\right),
\]
where $A = VH$ is the polar decomposition of $A$.
Since in the polar decomposition, we have $\lambda_i(H) = \sigma_i(A)$ where $\lambda_i(H)$ is the $i$th eigenvalue of $H$ and  $\sigma_i(A)$ is the $i$th singular value of $A$, the standard initialization provides an initial error of
\begin{equation} \label{boundstan}
{\| A - V \,p_{0,1}(H)\|}_F^2 = \sum_{i, \sigma_i(A) > 1} (\sigma_i(A)-1)^2.
\end{equation}

\paragraph{LMI-based initialization} Let
$\mu = \max(1,\rho(A))$ so that $A' = \frac{A}{\mu}$ is stable.
Then, we use the Lyapunov solution $P \succ 0$ to the system ${A'}^TPA'-P \preceq 0$ (we used the Matlab function \texttt{dlyap(A,eye(n)}) and define $S = P^{1/2}$, $R = S A' S^{-1}$, and $(U,B)$ is the polar decomposition of $R = UB$ so that $A' = S^{-1} U B S$; see the proof of
Theorem~\ref{thm:mainresult}.  This initialization provides a solution with initial error:
\begin{equation} \label{boundlmi}
{\| A - A'\|}_F^2 = {\| A - A/\mu\|}_F^2 = {\| A \|}_F^2 (1 - 1/\mu)^2.
\end{equation}

\begin{remark}[Comparing \eqref{boundstan}  and \eqref{boundlmi}] \label{errboundsinit}
None of the two solutions from \eqref{boundstan}  and \eqref{boundlmi} dominate the other one. It depends on the singular- and eigen-values of $A$.
For example, $A$ may be stable so that $|\rho(A)| < 1$ while $\sigma_{\max}(A) $ (the largest singular value of $A$) is greater than one in which case \eqref{boundlmi} provides an optimal solution (with error zero) while \eqref{boundstan} has a positive error.
On the other hand, if $A$ is symmetric so that $\sigma_{\max}(A) = \rho(A)$ and $A$ has a single singular value larger than 1, then the solution~\eqref{boundstan} has smaller error than~\eqref{boundlmi}. In fact,
\begin{eqnarray*}
{\| A \|}_F^2 (1 - 1/\mu)^2
&=&
\sum_{i=1}^n \left(\sigma_i(A)\right)^2 \left( \frac{\sigma_{\max}(A)-1}{\sigma_{\max}(A)} \right)^2 \\
&\geq& \left(\sigma_{\max}(A)\right)^2 ( \sigma_{\max}(A)-1 )^2\\
&>& ( \sigma_{\max}(A)-1 )^2.
\end{eqnarray*}
\end{remark}

\paragraph{Random initialization} We generate each entry of $S$ using the normal distribution (in Matlab, \texttt{randn(n)}).
Then, we replace $S$ with $SS^T+I_n$ which is positive definite.
Ideally, we then would like to compute the corresponding optimal $(U,B)$, that is, minimize ${\|A - S^{-1} UB S\|}_F$. However, we do not know how to do this efficiently, and instead we take $U$ and $B$ as the optimal solution of
\[
\min_{ U \text{ orthogonal}, B \succeq 0, \|B\| \leq 1} \; \|S A S^{-1} - UB\|_F,
\]
that is, $(U,B)$ is the polar decomposition of $S A S^{-1}$ and $B$ is replaced with $p_{0,1}(H)$; see Proposition~\ref{thm:optimUB}.
The motivation is that if $S A S^{-1} \approx  UB$ then $A \approx  S^{-1} UB S$. In fact, ${\|S X S^{-1}\|}_F \leq \kappa(S) {\|X\|}_F$ for any $X$. \\
In general, using a single random initialization provides a poor solution compared to the two previously proposed initializations. However, we have developed a simple multi-start heuristic that works as follows. Given a total allotted time $t_{\max}$ to the algorithm, we spend $t_{\max}/2$ generating 100 random initializations and refine them using Algorithm~\ref{fastgrad} (which therefore runs for only $t_{\max}/200$ for each random initialization). Then, we keep the best solution obtained among the 100 random initializations and refine it for $t_{\max}/2$.

\section{Numerical experiments}  \label{sec:Numerical}

In this section, we compare our algorithm, which we refer to as FGM,
with the only other known method for solving~\eqref{eq:probdef},
namely the successive convex approximation approach~\cite{OrbNV13}, kindly made available to us by Fran\c{c}ois-Xavier Orban de Xivry, that we refer to as SuccConv. \\
For both methods, we will use the standard and the LMI-based initializations.
FGM initialized with the standard (resp.\@ LMI-based) initialization is denoted Stand-FGM (resp.\@ LMI-FGM), and similarly for SuccConv.
We will use the multi-start heuristic only for FGM, which we will refer to as mRand-FGM,
because it is not well suited for SuccConv that required much more time per iteration, and more iterations to converge. \\
Our code is available from \url{https://sites.google.com/site/nicolasgillis/} and the numerical examples presented below can be directly run from this online code (there are also more numerical results in particular on randomly generated matrices). All tests are preformed using Matlab
R2015a on a laptop Intel CORE i5-3210M CPU @2.5GHz 6Go RAM.
FGM runs in $O(n^3)$ operations per iteration,
including projections on the set of positive semidefinite matrices, orthogonal matrices, and inversion of the matrix $S$ and all necessary matrix-matrix products.
Hence FGM can be applied on a standard laptop with $n$ up to a thousand. SuccConv is a second-order method and cannot be applied to matrices with $n$ much larger than 50
(one iteration of the algorithm requires about 30 seconds for $n=50$).

\subsection{Examples from~\cite{GP2018}}

We start with some examples from the paper~\cite{GP2018}.
In~\cite{GP2018}, authors study the nearest stable matrix problem~\eqref{eq:probdef} with component-wise nonnegativity constraint on the stable matrix $X$ to be found. For these small examples,
we set the time limit of the different algorithms to 30 seconds.

\subsubsection{Example 2: 3-by-3 matrix}

We consider
\[
A =\left( \begin{array}{ccc}
0.6 & 0.4&  0.1 \\
 0.5 & 0.5 & 0.3\\
 0.1 & 0.1 & 0.7
 \end{array} \right)
 \]
\normalsize  with $\rho(A) = 1.096$,  for which \cite{GP2018} shows that the nearest stable nonnegative matrix is
 \[
X =\left( \begin{array}{ccc}
0.5640 & 0.3599 & 0.0850 \\
 0.4716 & 0.4684 & 0.2881 \\
 0.0643 & 0.0602 & 0.6851
 \end{array} \right).
 \]
 \normalsize FGM and SuccConv for any initialization strategy converge to the same solution.
This is because, as shown in~\cite{GP2018} for nonnegative matrices, if a local minimum to problem~\eqref{eq:probdef}
is component wise positive, then it is a global minimizer.


\subsubsection{Example 3: scaled all-one matrix}

We now consider the matrix $A = \alpha E$ where $\alpha \geq 0$ and $E$ is the matrix whose entries are all equal to one. Note that the matrix is of rank-one, and stable for $\alpha \leq \frac{1}{n}$.
Authors~\cite{GP2018} show that for any $\alpha \in \left[ \frac{1}{n}, \frac{2}{n} \right]$, the nearest stable matrix is given by $\frac{1}{n} E$.
We run FGM and SuccConv on this example for $n = 10$ and $\alpha = \frac{2}{n}$ and,
as for the previous example,
they converge to the solution $\frac{1}{n} E$ for any of the initializations --note that LMI-FGM is initialized with the optimal solution since
$\frac{1}{n} E = \frac{A}{\rho(A)}$. (The same observation holds for $n=20$.)
For $\alpha > \frac{2}{n}$, $\frac{1}{n} E$ is not optimal anymore, and the optimal solution is not positive anymore~\cite{GP2018}.
For example, for $n = 2$ and any $\alpha > 1$, there are two nonnegative optimal solutions given by
$
\left( \begin{array}{cc}
1 & \alpha \\
0 & 1
 \end{array} \right)$,
 and
$\left( \begin{array}{cc}
1 & 0 \\
\alpha & 1
 \end{array} \right),
 $
 with error ${\|A-X\|}_F^2 = 6$.
Taking $\alpha = 3$,
LMI-FGM is not able to recover an optimal solution: it recovers $\frac{1}{n} E$ with error 9 since it is initialized with this solution and it is a stationary point of the problem~\cite{GP2018}.
Stand-FGM recovers a slightly better solution with error 8.
Stand-SuccConv and LMI-SuccConv obtain better but non-optimal solutions with error 6.27 and 6.24, respectively.
Only mRand-FGM is able to recover one of the above optimal solutions: rather surprisingly, it seems the unconstrained solution coincides with the nonnegative one (although we do not have a proof for this fact) --as we will see in the next examples, this is usually not the case.
For $n = 3$ and $\alpha = 2$,
the algorithms converge to different stationary points.
The triangular matrix with ones on the diagonal and $\alpha$ above or below has error 15.
As before, LMI-FGM converges to the stationary point $\frac{1}{n} E$ with error 25, and Stand-FGM to a better solution with error 17.
Stand-SuccConv and LMI-SuccConv converge to two rather different solutions with errors 15.2548 and 15.2558 respectively, while mRand-FGM provides the suboptimal solution
 \[
X =\left( \begin{array}{ccc}
0.9969 &   1.4010  &  0.7688 \\
    0.5544  &  0.9878  & -0.6507\\
    1.2476  &  2.6740  &  1.0112\\
 \end{array} \right)
 \]
 \normalsize with error 15.02. This illustrates the fact that~\eqref{eq:probdef} is a difficult problem with many local minimizers.

\subsubsection{Example in Section 4.4}
We consider
\[
A =
\left( \begin{array}{ccccc}
 0.7 &  0.2 &  0.1 &  0.5 &  1 \\
 0.3 &  0.6 &  0.2 &  0.8 &  0.3 \\
 0.5 &  0.7 &  0.9 &  1 &  0.5 \\
 0.1 &  0.1 &  0.3 &  0.8 &  0.3\\
 0.8 &  0.2 &  0.9 &  0.3 &  0.2 \\
\end{array} \right)
\]
\normalsize
with $\rho(A) = 2.4$. The nonnegative solution provided by the authors with their algorithm is
\[
X_+ =
\left( \begin{array}{ccccc}
 0.3796 &  0.1797 &  0 &  0.5 &  0.7343 \\
 0 &  0.5791 &  0.0069 &  0.8 &  0.0274 \\
 0.0580 &  0.6719 &  0.6403 &  1 &  0.1334 \\
 0 &  0 &  0 &  0.8 &  0 \\
 0.4204 &  0.1759 &  0.6770 &  0.3 &  0 \\
\end{array} \right)
\]
\normalsize with error ${\|A-X_+\|}_F^2 = 1.2181$ (which is not necessarily optimal).
Stand-SuccConv and LMI-SuccConv converge to the same solution
\[
\left( \begin{array}{ccccc}
 0.5999 &  0.1317 &  -0.0882 &  0.5337 &  0.8834 \\
 0.2582 &  0.5864 &  0.0967 &  0.8512 &  0.2089 \\
 0.4469 &  0.6904 &  0.8242 &  1.0419 &  0.4257 \\
 -0.0828 &  -0.1243 &  -0.2132 &  0.8209 &  0.0595 \\
 0.7076 &  0.1273 &  0.7126 &  0.3255 &  0.0923 \\
\end{array} \right)
\]
with error 0.5709.
Stand-FGM, LMI-FGM and mRand-FGM converge to
three different solutions with errors 0.6053, 0.5808 and 0.5759, respectively.

\subsection{Grcar matrices}

Grcar matrices of order $k$ are a banded Toeplitz matrix with its subdiagonal set to $-$1 and both its main and $k$ superdiagonals set to 1.
For example, when $n=5$ and $k = 3$, we have the following Grcar matrix
\[
\left( \begin{array}{ccccc}
 1  &  1  &  1  &  0  &  0  \\
 -1  &  1  &  1  &  1  &  0  \\
 0  &  -1  &  1  &  1  &  1  \\
 0  &  0  &  -1  &  1  &  1  \\
 0  &  0  &  0  &  -1  &  1  \\
\end{array} \right).
\]
Grcar matrices have all their eigenvalues outside the unit ball.
Notice that the nearest nonnegative stable matrix is given by $\max(A,0)$ with error ${\|A - \max(A,0)\|}_F^2 = n-1$.
Table~\ref{tab:grcar} reports the results for $k=3$ and $n = 5,10,20,50$ with time limit of $t_{\max}= 30, 60, 120, 300$ seconds, respectively.
\begin{center}
 \begin{table*}[h!]
 \resizebox{.965\textwidth}{!}{\begin{minipage}{\textwidth}
 \begin{center}
\caption{Comparison of the algorithms for Grcar matrices $A$ of order $3$:
final relative error in percent, that is,
$\frac{{\|A-X\|}_F}{{\|A\|}_F}$ where $X$ is the stable approximation of $A$,
and, in brackets, the number of iterations performed. The best solution is indicated in bold. The second column reports the relative error in percent of the nearest nonnegative stable matrix (that is, $\sqrt{n-1}/{\|A\|}_F$).  }  \label{tab:grcar}
{\small{ \begin{tabular}{|c|c|c|c|c|c|c|}
 \hline $n$ &  $\max(A,0)$   &  Stand-FGM & LMI-FGM & mRand-FGM & Stand-SuccConv & LMI-SuccConv    \\
 \hline
  5 &   47.14  &   \textbf{31.23} (5078) & \textbf{31.23} (5599) & 31.24 (13029) & 31.63 (9092)  & 31.64 (9104)  \\
 10 &   45.75  &   \textbf{30.02} (112539) & 33.08 (115262) & 33.18 (55136) & 30.88 (5188)  & 31.33 (5163)  \\
 20 &   45.20  &   41.64 (49225) & 45.34 (45417) & 46.51 (24539) & 40.07 (419)  & \textbf{39.41} (421)  \\
 50 &   44.91  &   \textbf{53.25} (34054) & 55.98 (35473) & 49.70 (16596) & 60.26 (6)  & 54.28 (6)  \\
\hline \end{tabular}}}
 \end{center}
  \end{minipage}}
 \end{table*}
 \end{center}
We observe that  Stand-FGM performs the best for $n=5,10,50$ and LMI-SuccConv for $n = 20$.
In most cases, the algorithms initialized with different initial points converge to different stationary points.

\begin{figure}[ht!]
\begin{center}
\begin{tabular}{ll}
\includegraphics[width=0.5\textwidth]{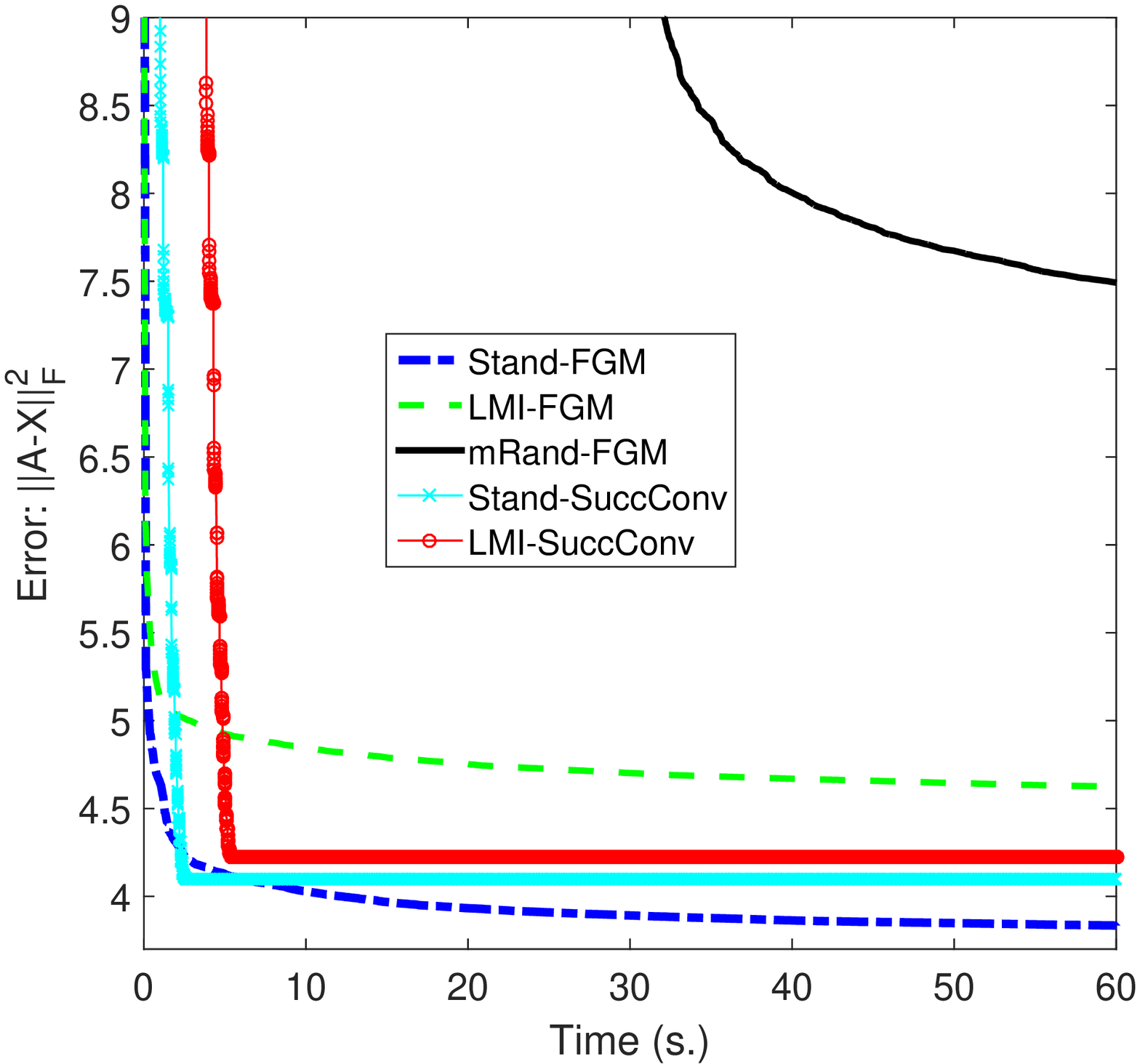} &
\includegraphics[width=0.5\textwidth]{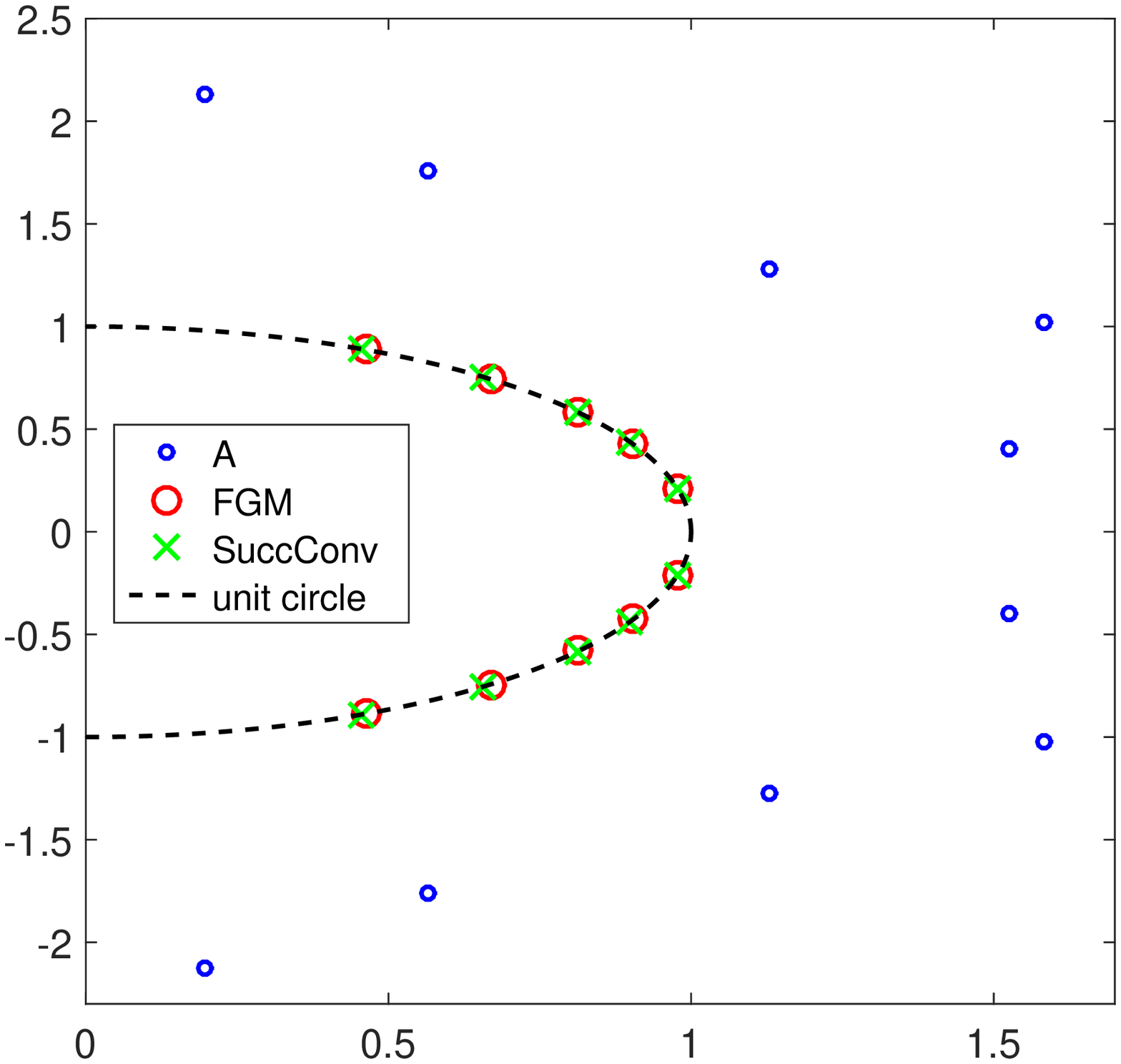}
\end{tabular}
\end{center}
 \caption{
(Left) Evolution of the error $\|A-X\|_F^2$ for the different algorithms for the Grcar matrix of dimension 10 and order 3. (Note that mRand-FGM only starts around 30 seconds as the multi-start heuristic spend half the time identifying the best solution among 100 randomly generated matrices.)
(Right) Location of the eigenvalues of $A$ and of the solutions obtained by FGM and SuccConv with the standard initialization.   \label{fig:grcar}}
\end{figure}

Figure~\ref{fig:grcar} shows the evolution of the objective function for the different algorithms for $n=10$ (on the left), and the location of the eigenvalues of $A$, and of the solutions of Stand-FGM and Stand-SuccConv, the best solution found by the two algorithms (on the right).
Although the eigenvalues
of the solutions $X_{\text{fgm}}$ of Stand-FGM and  $X_{\text{sc}}$ of Stand-SuccConv are close to one another, they actually correspond to very different matrices, since
$\frac{\|X_{\text{fgm}} - X_{\text{sc}}\|_F}{\|A\|_F} = 22.1\%$.

\section{Conclusion}

In this paper, we have proposed a new characterization of the set of stable matrices in the discrete-time case: We have shown that a matrix  $A$ is stable if and only if it admits a SUB form, that is, if there exists $S \succ 0$, $U$ orthogonal and $B \succeq 0$ with $\|B\| \leq 1$ such that
\[
A \; = \; S^{-1} U B S.
\]
We have then used this characterization to provide a new algorithmic framework for the nearest stable matrix problem, that is, given an unstable matrix $A$, find the nearest stable matrix $X$. In fact, the SUB form is particularly useful as it is easy to project onto this set of matrices.
We showed on several examples that our proposed algorithm that uses a fast gradient method (FGM) competes favorably with the method from~\cite{OrbNV13}.
In fact, in most cases, it provides better solutions while converging much faster.

Further research on the nearest stable matrix problem include the design of
(1)~other algorithms,
(2)~other initializations strategies, and
(3)~other heuristics to identify good solutions.
Further research also includes the use of the SUB form in defining the structure of linear port-Hamiltonian systems at the discrete level analogous to the continuous-time linear port-Hamiltonian systems, see, e.g.,~\cite{GolSBM03,Sch06,SchM13}, and to obtain the counterparts of the results in~\cite{MehMS16,MehMS17} for the discrete-time case.

\section*{Acknowledgments}

The authors would like to thank the reviewers for their insightful comments which helped improve the paper significantly.

\appendix
\section{Gradient with respect to $S$}  \label{appgradS}

The standard inner product on $\R^{n,n}$ is defined by $\langle A|B \rangle:={\rm tr}(A^T B)=\sum_{i,j} a_{ij}b_{ij}$.
The associated norm is the Frobenius norm,
$\|A\|_F= \sqrt{\langle A|A \rangle}=\sqrt{\sum_{i,j}a_{ij}^2}.$
The   relations $(AB)^T =B^T A^T$ and  ${\rm tr}(AB)={\rm tr}(BA)$ imply that
\begin{equation}\label{eq:leftrightmove}
\langle A|BC\rangle=\langle B^T A|C\rangle=\langle AC^T |B\rangle.
\end{equation}
Let ${\mathcal D}$ be a nonempty open subset of $\R^{n \times n}$.
A matrix $G\in \R^{n,n}$ is said to be the gradient of a function
${\mathcal D}\ni S \mapsto f(S) \in \R$
at $S_0\in {\mathcal D}$ if
 \begin{equation}\label{eq:gradientdef}
\frac{d}{dt}f(S(t))|_{t=0}=\langle G, \dot S(0)\rangle
\end{equation}
for all differentiable curves $\R \ni t\mapsto S(t)\in {\mathcal D}$ with $S(0)=S_0$ and derivative $\dot S(t)$.
It it easily seen that there is at most one matrix $G$ with this property. Notation: $G=\nabla f(S_0)$.
In the derivation below we omit the argument $t$ and the index 0. Furthermore we simply write
$\dot f$ for the left hand side of (\ref{eq:gradientdef}).

For fixed square matrices $A,C$ we are going to determine the gradient of the function
$$f(S)=\|A-S^{-1}CS\|_F^2=\langle R-A|R-A\rangle ,$$
where $R:=S^{-1}CS$.
The derivative of $R$ along a differentiable curve is
\begin{eqnarray*}
\dot R&=&S^{-1}C\dot S-(S^{-1} \dot S S^{-1})CS=R\,S^{-1} \dot S -S^{-1} \dot S\,R,
\end{eqnarray*}
where the left equation follows from  the product rule and fact that the derivative of the function
$S \mapsto S^{-1}$ along a differentiable curve is $-S^{-1}\dot S S^{-1}$
(this  is obtained by differentiating the relation $S^{-1}S=I$).
Now, the derivative of $f$ along a differentiable curve can be calculated as
\begin{eqnarray*}
\dot f&=& \langle \dot R|R-A\rangle + \langle R-A|\, \dot R\, \rangle \\
&=&2 \, \langle R-A\,|\, \dot R\, \rangle \\
&=& 2 \, \langle R-A \, | \,R\,S^{-1} \dot S -S^{-1} \dot S\,R\, \rangle \\
&=&2 \, \langle (R\,S^{-1})^T (R-A)  - S^{-T}(R-A)R^{T}\,|\,\dot S \, \rangle \\
&& \hspace{5.8cm}~ \text{(by (\ref{eq:leftrightmove}))}\\
&=& 2\, \langle S^{-T} [R^T (R-A)  - (R-A)R^{T}]\;|\,\dot S \, \rangle .
\end{eqnarray*}
Thus, the gradient of $f$ at $S$ is
\begin{equation*}
\nabla f(S)=2 \,S^{-T} [R^T (R-A)  - (R-A)R^{T}].
\end{equation*}

\newpage 

\small

\bibliographystyle{siam}
\bibliography{GilKS}

\end{document}